\documentclass[12pt]{article}
\begin{document}
\newcommand{\bq}{\begin{equation}}
\newcommand{\fq}{\end{equation}}
\newcommand{\IM}{{\bf i}}
\newcommand{\JM}{{\bf j}}
\newcommand{\MM}{{\bf m}}
\newcommand{\ccc}{C}
\newcommand{\aaa}{{\cal A}}
\newcommand{\qqq}{{\bf Q}}
\newcommand{\rrr}{{\bf R}}
\newcommand{\Z}{\bf Z}
\newcommand{\ZHS}{{\bf Z}HS}
\newcommand{\peak}{\check}
\newcommand{\rtb}{\begin{flushright}
   \begin{picture}(0,0)
    \put(0,5){\line(1,0){5}}
    \put(0,10){\line(1,0){5}}
    \put(0,5){\line(0,1){5}}
    \put(5,5){\line(0,1){5}}
   \end{picture}
  \end{flushright}}
\newcommand{\disjoint}{\cup}
\newtheorem{theorem}{Theorem}
\newtheorem{lemma}[theorem]{Lemma}
\newtheorem{conjecture}[theorem]{Conjecture}
\newtheorem{defn}[theorem]{Definition}

\begin{titlepage}
\hskip 1cm

\begin{flushright}
QMW-PH-96-18\\
q-alg/9608007\\
\vspace{12pt}
July 1996
\end{flushright}

\hskip 1.5cm
\vfill
\begin{center}
{\LARGE Ohtsuki's Invariants are of Finite Type}\vfill
{\large A. Kricker$^1$ and B. Spence$^2$\\ \vspace{6pt} }
\end{center}

\vfill
\begin{center}
{\bf Abstract}
\end{center}
{\small
Using a vanishing condition on certain combinations of components
of the Jones polynomial for algebraically split links we show that Ohtsuki's invariants of integral homology three spheres
are of finite type. We further show that the corresponding manifold weight
system is given by the expected Lie algebraic 
construction.
}

 \vfill \hrule width 3.cm
{\footnotesize
\noindent
 $^1$ Department of Physics, Queen Mary and Westfield College, London E1 4NS,
UK. On leave from the School of Mathematical Sciences, University of Melbourne, Parkville
3052 Australia. Email: A.Kricker@qmw.ac.uk\\
\noindent$^2$ 
Department of Physics, Queen Mary and Westfield College, London E1 4NS,
UK.
\\
Email:
B.Spence@qmw.ac.uk.
\\   }
\end{titlepage}


\section{Introduction}

Inspired by work in the theory of Vassiliev invariants of
knots,
the theory of  invariants of three-manifolds has developed rapidly of late. In particular, there has been much work done on
{\it finite-type} invariants of
integral homology three-spheres (\ZHS s) (see the revised version of \cite{L1} for a survey and references up to March 1996). Much of this work generalises to
rational homology three-spheres with little added complexity.

Ohtsuki \cite{O2} introduced finite type invariants of
\ZHS s (see also \cite{Ga1}). They are
characterised as follows. Consider a pair consisting of a \ZHS\ $M$ and an
embedded, algebraically-split,  
$\pm1$-framed link $L$ in $M$. Write $M_L$ for the \ZHS\ obtained via Dehn
surgery on $L$. A ZHS-invariant $\lambda$ is said to be {\it of finite type of order $\le n$} if
\bq
\sum_{L'\subseteq L} (-1)^{\#L'} \lambda(M_{L'}) = 0
\fq
when $\#L>n$, where $\#L$ is the number of components of the link $L$. It is
sufficient for this to be true for link embeddings in $S^3$.
(Such an invariant is said to be of {\it finite type of order $n$} if further 
it is not of finite type of order $m$ for any $m<n$.)
Ohtsuki has also suggested a candidate sequence of finite-type \ZHS\   invariants
$\{ \lambda_n \}$ \cite{O1} (see also \cite{LW}). He extracted the $\{ \lambda_n \}$ as a certain number-theoretic limit of the
Reshetikhin-Turaev $sl_2$ invariant \cite{RT}, when 
presented 
for \ZHS s. His work generalised H.Murakami's \cite{Mu} identification of
 Casson's invariant as the first non-trivial
invariant, $\lambda_1$.

In a separate important development, 
a universal finite-type \ZHS-invariant analogous to, and built from, Kontsevich's
construction for Vassiliev invariants \cite{K}, has been
presented in \cite{LMO}. It takes values in a certain vector
space built from trivalent graphs. It was established in \cite{Le'} that 
this invariant has a universal property for finite-type invariants of
\ZHS s (although it is defined for arbitrary manifolds).

An outstanding problem is to present the
$\{ \lambda_n \}$ invariants of Ohtsuki in 
terms of the universal invariant of \cite{LMO}. This would be analogous to 
the well-known
presentation of the quantum group knot invariants in terms of the Kontsevich
construction \cite{KT,LeMurakami}.

First, we will show here that the invariants $\{ \lambda_n \}$
are of finite type. Evidence for this has already been given
by Rozansky \cite{R1,R2},
using arguments from perturbative Chern-Simons theory. 
Our discussion utilises an explicit Kontsevich-type construction,
which we believe will also contribute 
towards a solution of the problem noted above. We will further confirm that
the manifold weight system of each Ohtsuki invariant is given in terms of
the $sl_2$ weight system for trivalent graphs. This provides substantial
support for Conjecture 7.3 of \cite{LMO}, namely that the 
series $\sum_{n=0}\lambda_n(e^h-1)^n$ is recovered from the 
universal finite-type
\ZHS\ invariant by ``substituting"  $sl_2$ into the graphs in the image of that
invariant.


\section{The Finite-Type Argument}

In the following we will consider invariants to be linearly extended
over formal linear combinations.
Namely, if we have some formal sum over links $C=\sum_{i=1}^n c_i L_i$,
and some link invariant $I$, then
$I(C)=\sum_{i=1}^n c_i I(L_i)$, with a similar construction for manifold invariants. We shall
also write  $M_C= \sum_{i=1}^n c_i M_{L_i}$.

The Ohtsuki invariants $\{ \lambda_n \}$ are constructed from derivatives of the Jones polynomial for framed knots. Define $V(L,h)$ to be the $\qqq[[h]]$ valued knot invariant defined by
\begin{eqnarray}
\label{Vdef}
e^h V(L_+,h) - e^{-h} V(L_-,h) &=& ( e^{\frac{h}{2}} - e^{-\frac{h}{2}} ) V(L_0,h), \nonumber \\
 V(\phi,h) &=& (e^{{h\over 2}} + e^{-{h\over 2}})^{-1},\ \    V(O,h)=1.
\end{eqnarray}
As usual, $L_+,L_-$ and $L_0$ are three links which differ only in a ball where
they are an overcrossing, an undercrossing, and an orientation-preserving 
smoothing of that crossing, respectively. $O$ is the unknot. The following definitions will prove useful (we follow the notation of \cite{LW} and \cite{L1}). Let
\begin{eqnarray}
X(L,h) & = & (e^{\frac{h}{2}}+e^{-\frac{h}{2}})^{1-\#L} V(L,h) \nonumber \\
& = & [2]^{1-\#L} V(L,h),
\end{eqnarray}
where the quantised integer $[n]$ is defined by
\bq
[n] = \frac{e^{hn/2} - e^{-hn/2}}{e^{h/2} - e^{-h/2}},
\fq
and $\#L$ is the number of components of the link $L$.
Now, for some link L, define
\begin{eqnarray}
\Phi(L,h) & = & \sum_{L' \subseteq L} (-1)^{\#L-\#L'} X(L',h) \\
& \equiv & \sum_{n=0}^{\infty} \frac{\Phi_n (L)}{n!} (e^h-1)^n.
\end{eqnarray}
Ohtsuki's invariants $\{ \lambda_n \}$ are constructed
using the functions $\phi_i(L)$ defined by
\bq
\phi_i(L) = \frac{(-2)^{\#L}}{(\#L+i)!} \Phi_{\#L+i}(L).
\fq

We will use the multi-index notation employed by Ohtsuki in \cite{O1},
which utilises bold-face letters for multi-indices, with the following
conventions. For ${\bf j}=(j_1,j_2,...,j_\mu)$ and 
${\bf k}=(k_1,k_2,...,k_\mu)$, and a scalar $x$, we put
\begin{eqnarray}
 |{\bf k}| &=& k_1+k_2+...+k_\mu, \nonumber\\
 \#{\bf k} &=& \mu, \nonumber \\
 {\bf j}+{\bf k } &=& (j_1+k_1,j_2+k_2,...,j_\mu+k_\mu), \nonumber\\
{\bf k}{\bf j} &=& (k_1j_1,k_2j_2,...,k_\mu j_\mu), \nonumber\\
 x{\bf k}  &=& (xk_1,xk_2,...,xk_\mu), \nonumber\\
{\bf x}  &=& (x,x,...,x).
\end{eqnarray}

The characterisation of the $\{\lambda_n\}$ we shall be exploiting here
is due to Lin and Wang \cite{LW}. In the phrasing of the theorem, we shall
denote by $L^n$ the zero-framed $n$th-parallel of the link $L$. That is,
we take the link that we get by replacing each component by $n$ parallel
components, all having linking number zero with each other, and each with the
same framing as the original component. Suppose that $L$ has the number of components $\#L=\mu$. We index 
sublinks  $L^{\bf i}\subseteq L^n$
by $\mu$-tuplets ${\bf i}=(i_1,\ldots,i_\mu)$, where $i_p$, for $p=1,...,\mu$, counts the number of
parallels of the $p$-th component of $L$ which appear in the sublink. 
(Note that in the sum over sublinks of $L^l$ below, each  
distinct $\mu$-tuplet is only counted {\it once}.) Then we have the
following formulation of Ohtsuki's invariants:
\begin{theorem}[\cite{LW}]
Let $S^3_L$ be the integral homology 3-sphere resulting from Dehn surgery on some
algebraically-split, $\pm1$-framed link $L$ embedded in $S^3$. Then there exist
constants $\{\nu_{f_p,i_p,m_p}\}$ such that
\bq
\label{Lin}
\lambda_n(S^3_L) = \sum_{l=1}^n \sum_{L'\subseteq L^l} \frac{\phi_l(L')}{(-2)^{\#L'}}
\left\{ \prod_{q=1}^\mu(-f_q) \left( \sum_{m_1+\ldots+m_\mu=n-l}\;
\prod_{p=1}^{\mu}\nu_{f_p,i_p,m_p}\right)\right\}.
\fq
\end{theorem}
\

\

In the above, $f_{q}=\pm1$ is the framing of the $q-$th component of $L$, 
and the
constants $\{\nu_{f_p,i_p,m_p}\}$ are defined recursively in \cite{LW}. Their
precise form in general will not be of direct concern to us here.

The following function will prove useful
\bq
F^{nl}(L,\IM,\JM) =  
\frac{1}{(-2)^{|\JM|}} \prod_{q=1}^{\mu} (-f_q)^{i_q}
\sum_{|\IM\MM|=n-l} \; \prod_{p=1}^{\mu} ( \nu_{f_p,j_p,m_p})^{i_p},
\fq
where in the summation $m_k$ is fixed at $0$ when $i_k=0$. 
Ohtsuki's invariants can now be expressed as
\bq
\lambda_n(S^3_{L}) = \sum_{l=1}^{n}\; \sum_{\JM={\bf 0}}^{l\bf 1}
\phi_{l}(L^{\JM})F^{nl}(L,{\bf 1},\JM).
\fq
We note the following technical identity
\begin{eqnarray}
\label{identityness}
F^{nl}(L,\{1,i_2,\ldots,i_\mu\},\{0,j_2,\ldots,j_\mu\}) &=& 
F^{nl}(L,\{0,i_2,\ldots,i_\mu\},\{0,j_2,\ldots,j_\mu\}). \nonumber\\
&&
\end{eqnarray}
This requires the result (see \cite{LW}) $\nu_{f,0,m}=-f\delta_{m0}$.
The identity (\ref{identityness}) follows from
\begin{eqnarray}
LHS\, (\ref{identityness}) & = & \frac{1}{(-2)^{|\JM|}} \prod_{q=2}^{\mu} (-f_q)^{i_q}(-f_1) \sum_{m_1+i_2m_2+\ldots+i_\mu m_\mu=n-l} \nu_{f_1,0,m_1}\prod_{p=2}^{\mu} (\nu_{f_p,j_p,m_p})^{i_p}\nonumber  \\
& = & \frac{1}{(-2)^{|\JM|}} \prod_{q=2}^{\mu} (-f_q)^{i_q} (-f_1)^2\sum_{m_2i_2+\ldots+m_{\mu}i_{\mu}=n-l} \;\;\prod_{p=2}^{\mu} (\nu_{f_p,j_p,m_p})^{i_p}\nonumber \\
 &=&  \frac{1}{(-2)^{|\JM|}} \prod_{q=2}^{\mu} (-f_q)^{i_q}
\;\; \sum_{m_2i_2+\ldots+m_\mu i_\mu=n-l} \;\;\prod_{p=2}^{\mu} (\nu_{f_p,j_p,m_p})^{i_p} \nonumber\\
    &=& RHS \,(\ref{identityness})
\end{eqnarray}

where we have used the fact that as all framings are $\pm 1$, they square to one. 

To prove that the invariants $\{\lambda_n\}$ are of finite type we need the following theorem
on the vanishing of the $\phi_i$.  
\begin{theorem}
\label{vanish}
Let $L$ be an algebraically split link. If $\#L > 3n$, then 
\bq
\phi_n(L)=0.
\fq
\end{theorem}
(This result for $n=1$  appeared in \cite{Hoste}, and the result
for $n=2$ was quoted in \cite{L1}.)
The proof of Theorem (\ref{vanish}) will occupy the next section.

With this understanding, we can turn to the proof of the main theorem.

\begin{theorem}
$\lambda_n$ is finite type of order $\leq 3n$.
\label{bigone}
\end{theorem}

{\bf \underline{Proof.}}

Suppose we have a link with $\mu=\#L>3n$. Then
\begin{eqnarray}
\sum_{L'\subseteq L} (-1)^{\# L'} \lambda_n(S^3_{L'}) &  = &
\sum_{\IM={\bf 0}}^{\bf 1} (-1)^{|\IM|}\;\sum_{l=1}^n \;\sum_{\JM={\bf 0}}^{l \IM} \phi_l(L^{\JM})\,
F^{nl}(L,\IM,\JM), \label{this}\\
& \equiv & \sum_{l=1}^n\; \sum_{\JM = {\bf 0}}^{l\bf 1} \phi_l(L^{\JM})\,G^{nl}(L,\JM), \label{thatonethere}
\end{eqnarray}
where (\ref{thatonethere}) {\it defines} the function $G$ as the coefficient of the appropriate $\phi$ in the  expression. Now we show that if at least one of $j_p$ are zero, then
$G^{nl}(L,\JM)=0$. Take some $\mu$-tuplet $\JM$ whose first $m$
indices are zero, with the remaining indices non-zero. From the
 above we have
\begin{eqnarray}
& & G^{nl}(L,\{0,\ldots,0,j_{m+1},\ldots,j_{\mu}\}) \nonumber \\
& = & \sum_{ \{i_1,\ldots,i_m \} = {\bf 0} }^{\bf 1} (-1)^{|\IM|} F^{nl}(L,
\{i_1,\ldots,i_m,1,\ldots,1\},\{0,\ldots,0,j_{m+1},\ldots,j_\mu\}),
\nonumber \\
&=&0,
\label{somany}
\end{eqnarray}
where in the final step we have used (\ref{identityness}). 
Hence
\begin{equation}
\sum_{L'\subseteq L} (-1)^{\# L'} \lambda_n(S^3_{L'})
 = (-1)^{\#L}\sum_{l=1}^n \sum_{\JM={\bf 1}}^{l\bf 1} \phi_l(L^{\JM})\,
F^{nl}(L,{\bf 1},\JM). \label{theresmore}
\end{equation}
We see then that (\ref{theresmore}) is a sum of terms, each with an expression
$\phi_l$ with $l\leq n$,  evaluated on links with more than $3n$ components. These vanish by Theorem 2.
Thus the expression on the left-hand side of (\ref{theresmore}) vanishes when $\#L>3n$, and we see that $\lambda_n$ is of finite type of order $\leq 3n$.
\rtb
\
  
\

This argument provides an upper bound for the order of
$\lambda_n$. Note from Theorem \ref{vanish} that when $\#L=3n$ the expression on the right-hand side of eqn.(\ref{theresmore})  has precisely one term, which is when ${\bf j}={\bf 1}$ and $l=n$ in the double
summand. Using the
fact that $\nu_{f,1,0}=2$ (see \cite{LW}), one then finds that for
the case when  $\#L=3n$,
\begin{eqnarray}
\label{umahhh}
\sum_{L' \subseteq L}(-1)^{\#L'} \lambda_n( S^3_{L'})
                 & = & 
(-1)^n\phi_n(L^{\bf 1})F^{nn}(L,{\bf 1},{\bf 1})  \nonumber \\
       &=&(-1)^nf_L\phi_n(L), 
\end{eqnarray}
where $f_L=\prod_{p=1}^{3n} f_p$.
We will use eqn.(\ref{umahhh}) in Section 4 to show that the
$\lambda_n$ are of order {\it exactly} $3n$.


\section{Proof of Theorem \ref{vanish}}

We begin by presenting the $\{ \phi_i \}$ in terms of Kontsevich's universal
invariant \cite{K}. 
Kontsevich's invariant is a chord diagram invariant, presented initially for 
knots but subsequently generalised to links and tangles. 
Different normalisations of the invariant exist, which are useful in different contexts. 
We will use  the tangle invariant $\hat{Z}$  as it appears in \cite{LMO,LMH} . This
normalisation factors the
quantum group tangle invariants of Reshetikhin and Turaev \cite{KT} . We denote by $\nu$ the value of $\hat{Z}$  on the unknot. 
We will also use the invariant $\hat{Z}'$, related to $\hat{Z}$ by
\bq
\hat{Z}(L)=(\nu)^{\mu}\#\hat{Z}'(L),
\fq
meaning we connect-sum a copy of
$\nu$ into each Wilson loop in the image of $\hat{Z}'(L)$.

Now, in the usual way,  choosing a Lie algebra with a
non-degenerate, symmetric, invariant metric, we can
define a weight system mapping the space of
chord diagrams $\hat{\cal A}$ into the rationals ${\bf Q}$, $\gamma: \hat{\cal A} \rightarrow {\bf Q}$. This is defined by
substituting the Lie algebra into the graph, and the representation into
the Wilson loop, and the grade of a chord diagram $D$ can be indexed in the 
weight system by defining $\Gamma: \hat{\cal A}\rightarrow {\bf Q}[[h]]$
by $\Gamma(D)=\gamma(D)h^{deg(D)}$  (see \cite{DBN} for example). In this work
we shall choose the algebra $sl_2$ with metric coming from the trace
${\rm Tr}_{fun}$ in the
{\it fundamental} representation $\rho$,
$<x,y>={\rm Tr}_{fun}(\rho(x)\rho(y))$, for $x,y\in sl_2$. 
We have 
\begin{lemma}[\cite{LMP}]
$\Gamma(\nu)=[2].$
\end{lemma}

The following result is a slight modification
of Theorem 2.3.1 of \cite{LMH}. It differs only in that here 
the invariant quadratic tensor which results from our choice of metric is 
one-half of
the tensor used to build weight systems in \cite{LMH}. This means that at grade
$n$, there is a multiplicative factor of $2^{-n}$. This simply scales the
appearance of $h$ in the statement as $h\rightarrow h/2$. Recall $V$ from eqn. (\ref{Vdef}).

\begin{theorem}[\cite{LMH}]
\bq
V(L) = \frac{1}{[2]}\Gamma(\hat{Z}(L_0)),
\fq
where $L_0$ is the zero-framed version of $L$, {\it i.e.} the link
$L$ with the framing of each of its
components set to zero.
\end{theorem}

We can use result this to obtain a formula for  $X$ as follows. (We use the simple identity $\Gamma(D_1\#D_2)=\Gamma(D_1)\Gamma(D_2)\frac{1}{2}$
for chord diagrams $D_1,D_2$.)
\begin{eqnarray}
\label{funnythat}
X(L) & = & [2]^{1-\mu} V(L) \nonumber\\
& = & \frac{1}{[2]^{\mu}} \Gamma(\hat{Z}(L_0)) \nonumber\\
& = & \frac{1}{[2]^{\mu}} \Gamma((\nu)^{\mu}\# \hat{Z}'(L_0)) \nonumber\\
& = & \frac{1}{2^\mu} \Gamma( \hat{Z}'(L_0) ).
\end{eqnarray}

Now we use this to obtain a formula for $\Phi$. For this we need two related operators.
In the following, a \lq leg\rq\ is an edge of the graph of some chord diagram that terminates in a
univalent vertex.

\begin{defn}
Choose a component C of L. 
Define an operator $\varepsilon_C: \hat{\cal A}\rightarrow \hat{\cal A}$
on chord diagrams $D$ by
\bq
{\varepsilon}_C(D) = \left\{ \begin{array}{ll}
0 & \mbox{if the graph of D has any legs on} \\
& \mbox{the Wilson loop corresponding to C,} \nonumber\\
D_C & \mbox{otherwise,}
\end{array} \right.
\fq
where $D_C$ is the diagram $D$ with the Wilson loop
corresponding to $C$ removed.

Define another operator $\tilde{\varepsilon}_C: 
\hat{\cal A} \rightarrow \hat{\cal A}$,  on chord diagrams by
\bq
\tilde{\varepsilon}_C(D) = \left\{ \begin{array}{ll}
D & \mbox{if the graph of D has any legs on} \\
& \mbox{the Wilson loop corresponding to C,} \nonumber \\
0 & \mbox{otherwise.}
\end{array} \right.
\fq
\end{defn}

\noindent
Then
\begin{lemma}[\cite{LMO}]
\bq
\hat{Z}'(L\backslash C)= \varepsilon_C(\hat{Z}'(L)),
\fq
\end{lemma}
and
\begin{lemma}
\begin{eqnarray}
X(L)-X(L\backslash C) & = &
\frac{1}{2^{\#L}}\Gamma(\hat{Z}'(L_0))  - \frac{1}{2^{\#L-1}}\Gamma(\hat{Z}'(L_0\backslash C)), \nonumber \\
& = &
\frac{1}{2^{\#L}}\Gamma(\hat{Z}'(L_0)  - 2\varepsilon_C(\hat{Z}'(L_0))), \nonumber \\
& = & \frac{1}{2^{\#L}}\Gamma(\tilde{\varepsilon}_C(\hat{Z}'(L_0))).
\end{eqnarray}
\end{lemma}
(We have used the fact that $\Gamma(O \disjoint D)=2\Gamma(D)$, where $O$ is the chord diagram with empty graph.)

Now notice that $\Phi(L)$ can be expressed as %
\bq
\Phi(L)= X(\delta_1\ldots\delta_\mu L),
\fq
 where $\delta_p L = L - L\backslash L_p$, with $L_p$ the $p$th
component of $L$.
Thus
\begin{eqnarray}
\Phi(L) & = &\frac{1}{2^\mu} \Gamma(\tilde{\varepsilon}_1(\ldots(\tilde{\varepsilon}_\mu
(\hat{Z}'(L_0))))), \nonumber \\
& = &\frac{1}{2^\mu} \Gamma(P(\hat{Z}'(L_0))),
\label{phi}
\end{eqnarray}
where $P$ is the projector which annihilates any chord
diagram which has a component without legs attached.

Finally, let us note that the Kontsevich invariant has special properties when evaluated on 
algebraically-split and zero-framed links (note that a {\it chord} is a 
component of the graph of a chord diagram without trivalent vertices):
\begin{lemma}
(i) Suppose that two components of some link $L$  have linking number zero. 
Then $\hat{Z}'(L)$ lies in the subspace generated by chord diagrams which
do not contain any
chords whose two ends terminate on the two Wilson loops corresponding to
the zero-linked components.\\
(ii) If some component of $L$ is zero-framed, then 
$\hat{Z}'(L)$  lies in the subspace generated by chord diagrams which do not
contain any chords
 whose two ends lie on the Wilson loop corresponding to the zero-framed
component.
\label{AS}
\end{lemma}
This Lemma is a straightforward  exercise in the theory of the Kontsevich
integral, and is implicit in various works in the literature (for example
see Section 5.4 in \cite{Le'}). We indicate the logic here. Every link can be
presented as the closure of a string link tangle (a tangle  where every 
component finishes vertically above its other end). The chord diagrams on
$n$-stringed string link tangles form an algebra generated by its 
primitive elements
({\it i.e.} those diagrams with connected graph), under the usual tangle 
multiplication. On string link tangles the Kontsevich invariant appears
as an exponential of primitive elements. It is straightforward to identify the 
coefficients of the diagrams with single chords in the logarithm of the
invariant with the appropriate linking numbers.

With this understanding, the proof of Theorem \ref{vanish} 
now follows from
a simple counting argument.

\

{\bf \underline{Proof of Theorem \ref{vanish}.}}

\

Consider a link $L$ with $\mu$ components. Recall that $\Phi$ is of the
form $\Phi(L)= \frac{1}{2^\mu} \Gamma(P(\hat{Z}'(L_0)))$. Observe from 
Lemma \ref{AS} that as the whole link $L_0$ is zero-framed and algebraically-split
then $\hat{Z}'(L_0)$ takes values in the subspace generated by chord diagrams
{\it without chords}. We then ask: what is the least grade, in this subspace, of a chord diagram
which is not annihilated by P, {\it i.e.} a diagram with legs on all 
components?

Every Wilson loop thus needs at least one leg terminating on it. Since there
are no chords, the other end
of this edge must end in a trivalent vertex. This implies that there
 must be at least {\it two}
legs ending on each Wilson loop, because of the anti-symmetry of that
trivalent vertex.

The least number of trivalent vertices that can provide two legs for
every component is $\frac{2\mu}{3}$. Such a diagram has grade $\frac{\mu}{3}+\mu$. Thus $\Phi$ evaluated on a link $L$ with $\mu$ components is of order at least
$h^{\mu+\frac{\mu}{3}}$. In other words
\bq
\phi_i(L)=0\ \ \mbox{if}\ 3i<\mu,
\fq
proving Theorem 2.
\rtb


\section{ The weight system of the invariants $\{ \lambda_n \}$}

In this section we shall show that the manifold weight system (or 
\lq corresponding linear form\rq\ \cite{Le'}) of $\lambda_n$ is 
evaluated by \lq substituting\rq\ $sl_2$ into certain trivalent graphs.
A corollary of this result is that the invariants $\{ \lambda_n \}$
are of order $3n$ exactly.

We denote by
${\cal M}$ the $\qqq$-vector space spanned by homeomorphism classes
of integral homology three-spheres. Consider an algebraically split,
$\pm 1$-framed link embedded in some \ZHS\ $M$. Define
\bq
\delta(L)=\sum_{L'\subseteq L} (-1)^{\# L'} L'.
\fq
Denote by ${\cal F}_n$ the subspace of ${\cal M}$ generated by vectors of
the form $M_{\delta(L)}$ where $\# L=n$. This provides a filtration on ${\cal M}$.
In this notation the condition 
that an invariant $I$ be of finite-type of order $\leq n$ is that $I({\cal F}_{n+1})=0.$

Let ${\cal D}_n$ be the $\qqq$-vector space spanned by trivalent 
vertex-oriented graphs
with $2n$ vertices, modulo the usual AS and IHX relations \cite{LMO}. In \cite{GO3} a surjective
linear mapping $O^*_{3n}: {\cal D}_n \rightarrow {\cal F}_{3n}/{\cal F}_{3n+1}$
was constructed (it was shown in \cite{Le'} that this mapping is an 
isomorphism). Thus invariants $I$ of order $3n$ give rise to well-defined manifold
weight systems via
\bq
I\circ O^*_{3n}: {\cal D}_n \rightarrow \qqq.
\fq
Here we shall consider the manifold weight systems arising from
the Ohtsuki invariants, given by 
\bq
\Lambda_n = \lambda_n \circ O^*_{3n}.
\fq
Let us first recall the map $O^*_{3n}$. Let $G\in {\cal D}_n$ be a trivalent
 vertex-oriented graph with $2n$ vertices and $3n$ edges. 
Present $G$ as a graph diagram
in the plane ({\it i.e.} embedded with  double points labelled with
crossing information, and cyclic orientation at trivalent vertices, assumed 
counter-clockwise). We construct a formal sum of $3n$-component links
from this (\cite{Le'}, \cite{GO3}). First, replace every vertex with a difference of
three-component links as in Fig. 1, with each of the link components being
associated to a different 
incoming vertex edge.

\begin{figure}[htp]
\centerline{
\begin{picture}(120,60)
\put(17,29){\line(1,0){5}}
\put(24,29){\line(1,0){5}}
\put(31,29){\line(1,0){5}}
\put(38,29){\line(1,0){5}}
\put(45,29){\line(1,0){5}}
\put(52,29){\line(1,0){5}}
\put(60,30){\circle*{2}}
\put(63,29){\line(1,0){5}}
\put(70,29){\line(1,0){5}}
\put(77,29){\line(1,0){5}}
\put(84,29){\line(1,0){5}}
\put(91,29){\line(1,0){5}}
\put(98,29){\line(1,0){5}}
\put(60,23){\line(0,1){5}}
\put(60,16){\line(0,1){5}}
\put(60,9){\line(0,1){5}}
\put(60,2){\line(0,1){5}}
\put(60,-5){\line(0,1){5}}
\put(60,-12){\line(0,1){5}}
\end{picture}
\hskip0cm
\begin{picture}(10,60)
\put(5,20){\makebox(0,0){$\longrightarrow$}}
\end{picture}
\hskip 0cm
\begin{picture}(120,60)
\put(33,51){\line(1,0){14}}
\put(42,36){\oval(72,30)[rt]}
\put(58,28){\oval(40,30)[rb]}
\put(33,13){\line(1,0){12}}
\put(52,13){\line(1,0){9}}
\put(40,36){\oval(40,30)[lt]}
\put(20,27){\line(0,1){13}}
\put(40,28){\oval(40,30)[lb]}
\put(63,26){\oval(30,14)[tl]}
\put(48,-6){\line(0,1){33}}
\put(71,26){\oval(30,14)[tr]}
\put(86,-6){\line(0,1){26}}
\put(63,-5){\oval(30,14)[bl]}
\put(71,-5){\oval(30,14)[br]}
\put(60,-12){\line(1,0){12}}
\put(71,36){\oval(8,20)[lt]}
\put(79,47){\line(1,0){29}}
\put(67,29){\line(0,1){10}}
\put(72,32){\oval(10,20)[lb]}
\put(81,22){\line(1,0){27}}
\put(107,37){\oval(8,20)[rt]}
\put(107,32){\oval(8,20)[rb]}
\put(111,29){\line(0,1){8}}
\end{picture}
\hskip 0cm
\begin{picture}(10,60)
\put(5,20){\makebox(0,0){$-$}}
\end{picture}
\hskip 0cm
\begin{picture}(120,60)
\put(31,37){\oval(64,30)[rt]}
\put(63,27){\line(0,1){14}}
\put(43,28){\oval(40,30)[rb]}
\put(33,13){\line(1,0){12}}
\put(40,37){\oval(40,30)[lt]}
\put(20,27){\line(0,1){14}}
\put(40,28){\oval(40,30)[lb]}
\put(63,0){\oval(30,14)[tl]}
\put(48,-7){\line(0,1){7}}
\put(61,7){\line(1,0){12}}
\put(71,0){\oval(30,14)[tr]}
\put(86,-7){\line(0,1){7}}
\put(63,-5){\oval(30,14)[bl]}
\put(71,-5){\oval(30,14)[br]}
\put(60,-12){\line(1,0){12}}
\put(76,36){\oval(8,20)[lt]}
\put(75,46){\line(1,0){28}}
\put(72,29){\line(0,1){10}}
\put(77,32){\oval(10,20)[lb]}
\put(76,22){\line(1,0){27}}
\put(103,36){\oval(8,20)[rt]}
\put(103,32){\oval(8,20)[rb]}
\put(107,29){\line(0,1){8}}
\end{picture}
}
\vskip 0.5cm\caption{The Borromean difference}
\end{figure}
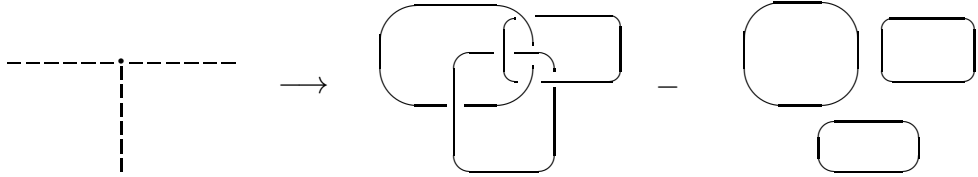

Then we connect sum pairs of link components whenever their
corresponding vertex edges connect two
trivalent vertices, with the joining pairs of arcs lying along the connecting edge. From
a graph with $2n$ vertices, this yields a sum of $2^{2n}$ links,
each with $3n$-components, taken to be zero-framed. Call this sum $\beta(G)$. Now alter the framing of each component to $+1$,
and call 
this formal sum of links $\tilde{\beta}(G)$. 
It is a simple exercise to show that $\delta(\tilde{\beta}(G)) = 
(-1)^n \tilde{\beta}(G)$. 
The vector $O^*_{3n}(G)$ is defined as the image of 
$S^3_{\tilde{\beta}(G)}$ in ${\cal F}_{3n}/{\cal F}_{3n+1}$.

The evaluation we require is 
\begin{eqnarray}
\Lambda_n(G) & =& \lambda_n(S^3_{\tilde{\beta}(G)}) \nonumber \\
& = & \lambda_n(S^3_{\delta((-1)^n\tilde{\beta}(G))}).
\label{blahblah}
\end{eqnarray}
{}From eqn. (\ref{umahhh}), noting that all framings are 1, we see that
\begin{eqnarray}
\Lambda_n(G) & = &  (-1)^n\phi_n((-1)^n\tilde{\beta}(G) ) \nonumber \\
& = &  (-2)^{3n} \frac{\Phi_{4n}( \tilde{\beta}(G))}{(4n)!} \nonumber \\
& = & (-1)^n Grad_{4n} ( \Gamma( P( \hat{Z}' ( \beta(G))))),
\label{answer}
\end{eqnarray}
where $Grad_{4n}$ means take the coefficient of $h^{4n}$, and we have the used
the equality (\ref{phi}).
In the last equation above, we have used the fact that here $\Phi_i(\tilde\beta(G))=0$ for $i<4n$ by Theorem 2.
Thus we can identify $\Phi_{4n}$ with 
the coefficient
of $h^{4n}$ (this is not the case in general, due to the different choices of expansion parameters).

An analysis of the series $\hat{Z} ( {\beta}(G) )$
has been performed in \cite{Le'}. Following this reference, we will say that a chord diagram has 
{\it i-filter $\geq$ n} if it is in the subspace spanned by chord diagrams with at
least $n$ trivalent vertices on their associated graphs. We will also
use a linear mapping
$\eta: {\cal D}_n \rightarrow {\cal A}(\coprod (S^1)^{3n})_{4n}$ defined simply
by breaking each edge of some graph $G\in {\cal D}_n$ and inserting a Wilson
loop, as in Fig. 2.

\begin{figure}[htp]
\centerline{
\begin{picture}(120,60)
\put(17,29){\line(1,0){5}}
\put(24,29){\line(1,0){5}}
\put(31,29){\line(1,0){5}}
\put(38,29){\line(1,0){5}}
\put(45,29){\line(1,0){5}}
\put(52,29){\line(1,0){5}}
\put(59,29){\line(1,0){5}}
\put(66,29){\line(1,0){5}}
\put(73,29){\line(1,0){5}}
\put(80,29){\line(1,0){5}}
\put(87,29){\line(1,0){5}}
\put(94,29){\line(1,0){5}}
\end{picture}
\hskip0cm
\begin{picture}(10,60)
\put(5,28){\makebox(0,0){$\rightarrow$}}
\end{picture}
\hskip 0cm
\begin{picture}(120,60)
\put(17,29){\line(1,0){5}}
\put(24,29){\line(1,0){5}}
\put(31,29){\line(1,0){5}}
\put(38,29){\line(1,0){5}}
\put(45,29){\line(1,0){5}}
\put(52,29){\line(1,0){5}}
\put(70,30){\circle{25}}
\put(84,29){\line(1,0){5}}
\put(91,29){\line(1,0){5}}
\put(98,29){\line(1,0){5}}
\put(105,29){\line(1,0){5}}
\put(112,29){\line(1,0){5}}
\put(119,29){\line(1,0){5}}
\end{picture}
}
\vskip 0.5cm\caption{Inserting a Wilson loop}
\end{figure}
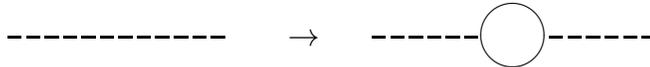

\begin{lemma}[\cite{Le'}, eqn. (5.4)]
\bq
\hat{Z}( \beta( G)) = \eta(G) + (\mbox{elements of i-filter}> 2n).
\fq
\label{graph}
\end{lemma}
The function $\hat{Z}'$ is obtained from $\hat{Z}$ by connect-summing a
$\nu^{-1}$ into each component.
Now, $\nu$ has the form $1 + ($elements of grade$\geq 1$), and hence its inverse will take the same form. Thus, at order $4n$,
Lemma \ref{graph} holds true if we replace $\hat{Z}$ by
$\hat{Z}'$. Similarily, choosing zero-framing for each component does not affect
the result (\ref{graph}) at grade $4n$ either, as this just multiplies $\hat{Z}'$ by some exponential 
of the chord diagram with one isolated chord.

Finally, note that at order $4n$, elements of i-filter$> 2n$ are in the
kernel of the projector $P$, as they have less than $6n$ legs (the diagram
requires $6n$ legs to put two legs on every Wilson loop).

Thus,
\bq
\label{doitnow}
\Lambda_n(G) = (-1)^n\Gamma( \eta( G) ).
\fq

Given that $\Gamma$ is defined from $sl_2$ with metric coming from the trace
in the fundamental representation, $\Gamma( \eta(G))$ is just equivalent
to substituting the Lie algebra into the graph $G$, in the usual fashion. 

Finally, noting that the  Lie algebraic weight systems are non-trivial,
we can see from Theorem 3 and equation (\ref{doitnow})
that the following is true:

\begin{theorem}
$\lambda_n$ is of finite type of order exactly $3n$.
\end{theorem}

At this point it is worth comparing  the weight system we have calculated
for the $\{ \lambda_n \}$ with what we would expect if Conjecture (7.3) of
\cite{LMO} were correct. Denote by  $\hat\Omega$ the universal invariant
of \cite{LMO}.
A specific case of their conjecture is
\begin{conjecture}[\cite{LMO}]
Let $M$ be an integral homology three sphere. Then
\bq
1+\sum_{m=1} \lambda_m (e^h - 1)^m = \Gamma(\hat\Omega(M)).
\fq
\end{conjecture}
The calculation of the weight system for the
universal invariant appears in \cite{Le'} (Lemma 5.2):
\begin{lemma}[\cite{Le'}]
For graphs $G\in {\cal D}_n$,
\bq
\Omega_n(S^3_{\tilde{\beta}(G)}) = (-1)^n G.
\fq
\end{lemma}
Thus our result, eqn. (\ref{doitnow}), represents significant support for 
the conjecture in the case $sl_2$. What one can ask further, by analogy
with the case for knots, is the question: Is $Grad_n(1+\sum_{m=0} \lambda_m (e^h-1)^m)$
canonical? 

We finish by noting that with the understanding of the $\{\phi_i\}$ 
presented in this work, the identification of the coefficient of the 
theta-curve in $\Omega_1$ with Casson's invariant is straightforward,
requiring only the identification of the various contributing summands
with known quantities (see \cite{Cas} for a more technical approach).

\

\

{\bf \underline{Acknowledgements.}}

AK acknowledges the support of an Australian Postgraduate Research Award, an International Collaborative Research Scholarship and support from the School of Mathematical Sciences, University of Melbourne. 
He also thanks the string theory group at QMW for their hospitality during the course of this work. BS was supported by
a UK EPSRC Advanced Fellowship.

\

\

{\bf Note Added}: After circulating this preprint, we were informed by
T. Ohtsuki that he has recently found a proof of Conjecture (12) in the setting of rational homology three-spheres.



\end{document}